%% file: paper-rss.tex
\documentclass[conference]{IEEEtran}
\usepackage{times}

\usepackage[numbers]{natbib}
\usepackage{multicol}
\usepackage[bookmarks=true]{hyperref}

\usepackage{amsmath}
\usepackage{amssymb}
\usepackage{mathtools}
\usepackage{xcolor}
\usepackage{float}
\usepackage{mathrsfs}
\usepackage{relsize}
\usepackage{url}
\usepackage{cleveref}
\usepackage{booktabs}
\usepackage{wrapfig}
\input{mynotation.tex}

\usepackage[algo2e, ruled, vlined, linesnumbered]{algorithm2e}

\SetKwFor{While}{while}{}{}%
\SetKwFor{GotoLoop}{}{}{}%
\SetKwFor{For}{for}{}{}%
\SetKwIF{If}{ElseIf}{Else}{if}{}{else if}{else}{end if}%
\newcommand\mykwcomment{\quad$\triangleright$\ }
\SetKwComment{Comment}{\mykwcomment}{}
\SetKwComment{FullLineComment}{$\triangleright$\ }{}

\SetCommentSty{mycommentstyle}

\begin{document}

\title{Exploiting Parallelism in a QPALM-based Solver for Optimal Control}

\author{\authorblockN{Pieter Pas\authorrefmark{1},\quad
Kristoffer L\o{}wenstein\authorrefmark{2}\authorrefmark{3},\quad
Daniele Bernardini\authorrefmark{3} \;\,and\,\;
Panagiotis Patrinos\authorrefmark{1}}
\authorblockA{\authorrefmark{1}Department of Electrical Engineering (ESAT), KU Leuven, Kasteelpark Arenberg 10, 3001 Leuven, Belgium.\\
Email: \texttt{\small\{pieter.pas,panos.patrinos\}@esat.kuleuven.be}}
\authorblockA{\authorrefmark{2}Dipartimento di Elettronica, Informazione e Bioingegneria, Politecnico di Milano, Via Ponzio 34/5, 20133, Milano, Italy.\\
Email: \texttt{\small kristofferfink.loewenstein@polimi.it}}
\authorblockA{\authorrefmark{3}ODYS S.r.l., Via Pastrengo 14, 20159, Milano, Italy.\\
Email: \texttt{\small\{kristoffer.lowenstein,daniele.bernardini\}@odys.it}}\\[-0.1em]
July 15, 2024\\[-0.5em]}

\maketitle

\begin{abstract}
We discuss the opportunities for parallelization in the recently proposed QPALM-OCP
algorithm, a solver tailored to quadratic programs arising in optimal control.
A significant part of the computational work can be carried out
independently for the different stages in the optimal control problem. We
exploit this specific structure to apply parallelization and vectorization
techniques in an optimized \Cpp{} implementation of the method. Results for
optimal control benchmark problems and comparisons to the original QPALM method are provided.
\end{abstract}

\IEEEpeerreviewmaketitle

\section{Introduction}
The need for the solution of linear-quadratic optimal control problems (OCPs) arises in
popular applications such as linear model predictive control (MPC) and moving
horizon estimation (MHE), as well as in sub-problems of solvers for their nonlinear counterparts. The real-time nature of these applications requires
efficient quadratic programming (QP) solvers, that often need to perform well
in constrained embedded environments.
This work considers the parallelization of an optimized implementation of the
recently proposed QPALM-OCP algorithm \cite{10550920}, which is a specialization
of the augmented Lagrangian-based QPALM solver for (possibly nonconvex)
quadratic programs \cite{qpalm}.
The contribution in \cite{10550920} is twofold: first, the QPALM algorithm is
modified to handle linear equality constraints directly, and second, the specific
structure of linear-quadratic optimal control problems (OCPs) is exploited in
the solution of the linear systems in its inner semismooth Newton solver.
Handling equality constraints directly reduces the overall number of iterations,
and combined with more efficient linear solves in the inner solver, this can
significantly reduce the solver run time.

This paper takes this one step further: it further explores the stage-wise
structure of OCPs, and uses this to efficiently parallelize the algorithms
described in \cite{10550920}. %
Two levels of parallelism are used to take full advantage of the performance of modern hardware: (1) a compact storage format of the OCP matrices
allows for efficient single instruction, multiple data (SIMD) operations
across stages, and (2) OpenMP \cite{openmp} is used to
distribute the stages across multiple physical CPUs.

The remainder of this document is structured as follows: \Cref{sec:formulation}
introduces the optimal control problem formulation and notation; \Cref{sec:qpalm-cop} summarizes the QPALM-OCP algorithm and highlights the opportunities for parallelization. \Cref{sec:vec,sec:par} discuss data parallelism using vectorization and OpenMP parallelization.
\Cref{sec:results} presents benchmark results, and \Cref{sec:conclusion} concludes the paper with discussions of future work.

\section{Problem formulation}\label{sec:formulation}

We consider linear-quadratic optimal control problems:
\begin{equation}
     \begin{aligned}
        &\minimize_{\mathbf{u}, \mathbf{x}}&& \sum_{j=0}^{N-1} \ell_j(x^j, u^j) + \ell_N(x^N) \\
        &\subjto && x^0 = x_\text{init} \\
        &&& x^{j+1} = A_j x^j + B_j u^j + c^j &&{\scriptstyle(0 \le j \lt N)} \\
        &&& b_l^j \le C_j x^j + D_j u^j \le b_u^j &&{\scriptstyle(0 \le j \lt N)} \\
        &&& b_l^N \le C_N x^N \le b_u^N
    \end{aligned} 
\end{equation}
The convex stage-wise and terminal cost functions are given by $\ell_j(x, u) = \tfrac12 \begin{pmatrix}
    x \\ u
\end{pmatrix}^{\!\top}\!\! \begin{pmatrix}
    Q_j & \tp S_j \\ S_j & R_j
\end{pmatrix} \begin{pmatrix}
    x \\ u
\end{pmatrix} + \begin{pmatrix}
    x \\ u
\end{pmatrix}^{\!\top}\!\! \begin{pmatrix}
    q^j \\ r^j
\end{pmatrix}$ and $\ell_N(x) = \tfrac12 \tp x Q_N x + \tp x q^N$, with $Q_j \in \possdefset{\R^{n_x}}$, $S_j \in \R^{n_u\times n_x}$ and $R_j \in \possdefset{\R^{n_u}}$.
The state variables $\mathbf{x} = [x^0\, \cdots\, x^N] \in \R^{n_x\times(N+1)}$ and inputs $\mathbf{u} = [u^0\, \cdots\, u^{N-1}] \in \R^{n_u\times N}$
are governed by the linear dynamics described by $A_j \in \R^{n_x\times n_x}$,
$B_j \in \R^{n_x\times n_u}$ and $c^j \in \R^{n_x}$. Linear stage-wise mixed state-input
constraints are encoded using $C_j\in\R^{n_y \times n_x}$ and $D_j\in\R^{n_y \times n_u}$,
with $C_N \in\R^{n_{y^N} \times n_x}$ for the terminal constraints.

By interleaving the states and inputs into a vector $x\in\R^{N(n_x + n_u) + n_x}$ and
introducing the matrices $H$, $M$ and $G$, the OCP can be written as a quadratic
program in standard form, with $n \defeq N(n_x + n_u) + n_x$ variables, $p \defeq (N+1) n_x$
equality constraints, and $m \defeq N n_y + n_{y^N}$ inequality constraints:%
\begin{equation} \label{eq:standard-qp}
    \begin{aligned}
        &\minimize_x && \tfrac12 \tp x Q x + \tp x q \\
        &\subjto && M x = b \\
        &&& b_l \le G x \le b_u \\
    \end{aligned}
\end{equation}
\def\shortminus{\scalebox{0.6}[1.0]{$-$}}
\begin{align*}
    Q &\defeq \scalebox{0.87}{$\setlength{\arraycolsep}{3pt}\begin{pmatrix}
        Q_0 & \tp S_0 \\
        S_0 & R_0 \\
        && Q_1 & \tp S_1 \\
        && S_1 & R_1 \\
        &&&& \ddots \\
        &&&&& Q_N
    \end{pmatrix}$}, \; x \defeq \scalebox{0.87}{$\setlength{\arraycolsep}{3pt}\begin{pmatrix}
        x^0 \\ u^0 \\ x^1 \\ u^1 \\ \vdots \\ x^N\hspace{-1pt}
    \end{pmatrix}$}, \;\;\, q \defeq \hspace{1pt}\scalebox{0.87}{$\setlength{\arraycolsep}{3pt}\begin{pmatrix}
        q^0 \\ r^0 \\ q^1 \\ r^1 \\ \vdots \\ q^N
    \end{pmatrix}$} \\
    M &= \scalebox{0.87}{$\setlength{\arraycolsep}{3pt}\begin{pmatrix}
        \I \\
        \shortminus A_0 & \shortminus B_0 & \I \\
        && \shortminus A_1 & \shortminus B_1 & \I \\
        &&&&& \ddots \\
        &&&& \shortminus A_{N-1} & \shortminus B_{N-1} & \I\,
    \end{pmatrix}$}, \;\; b \defeq\hspace{-1pt} \scalebox{0.87}{$\setlength{\arraycolsep}{3pt}\begin{pmatrix}
        x_\text{init} \\ c^0 \\ c^1 \\ \vdots \\ c^{N-1}\hspace{-1pt}
    \end{pmatrix}$} \\
    G &= \scalebox{0.87}{$\setlength{\arraycolsep}{3pt}\begin{pmatrix}
        C_0\;\;D_0 \\
        & C_1\;\;D_1 \\
        && \ddots \\
        &&& C_N
    \end{pmatrix}$}, \;\; b_{l} = \scalebox{0.87}{$\setlength{\arraycolsep}{3pt}\begin{pmatrix}
        b_{l}^0 \\
        b_{l}^1 \\
        \vdots \\
        b_{l}^N \\
    \end{pmatrix}$}, \;\;\! b_{u} = \scalebox{0.87}{$\setlength{\arraycolsep}{3pt}\begin{pmatrix}
        b_{u}^0 \\
        b_{u}^1 \\
        \vdots \\
        b_{u}^N \\
    \end{pmatrix}$}
\end{align*}

\section{QPALM-OCP}\label{sec:qpalm-cop}

The QPALM-OCP algorithm is described by \cite[Sec.~III]{10550920}, and further details
can be found in the main QPALM publication \cite{qpalm}. This section briefly
summarizes the inner semismooth Newton solver, which is responsible for the
main computational cost of the algorithm, and which will be optimized in the
remainder of this paper.

\subsection{The augmented Lagrangian inner problem}

QPALM-OCP solves problem \eqref{eq:standard-qp} by relaxing the inequality
constraints $b_l \le G x \le b_u$ using an augmented Lagrangian method (ALM).
The inner problem at outer iteration $k$ of this method amounts to the minimization
of a piecewise quadratic augmented Lagrangian function, subject to
equality constraints,
\begin{equation}\label{eq:inner-qp}
    \begin{aligned}
        &\minimize_{x} && \phi_k(x) \\[-0.2em]
        & \subjto &&Mx = b,
    \end{aligned}
\end{equation}
where $\phi_k(x) \defeq \tfrac12 \tp x Q x + \tp x q + \tfrac12 \dist^2_{\Sigma_{y,k}}\big(Gx + \inv \Sigma_{y,k} y^k, \, C\big) + \tfrac12 \normsq{x - x^k}_{\inv\Sigma_{x,k}}$.
In this expression, $y^k \in \R^m$ represents the current estimate of the Lagrange multipliers corresponding to the inequality constraints, $x^k \in \R^n$
is the current estimate of the primal variables, $\Sigma_{y,k} \in \posdefset{\R^m}$ is a diagonal matrix containing the ALM penalty weights,
the set $C = \defset{z\in \R^m}{b_l \le z \le b_u}$ describes the feasible set for $Gx$, and $\Sigma_{x,k} \in \posdefset{\R^n}$ adds primal regularization through a proximal term.

Because of the squared distance, the function $\phi_k$ is not twice
continuously differentiable. Therefore, QPALM applies a semismooth Newton method,
with a generalized Hessian matrix of $\phi_k$ that is given by \cite[Eq.~15]{10550920}
\def\setJ{\mathcal{J}}
\begin{equation}
    H_k(x) = Q + \tp G \Sigma_{y,k}^{\setJ\!(x)} G,
\end{equation}
where we defined the diagonal matrix $\Sigma_{y,k}^{\setJ\!(x)}$ as
\begin{equation}
    \big(\Sigma_{y,k}^{\setJ\!(x)}\big)_{ii} \defeq \begin{cases}
        0 & \text{if } \big(Gx + \inv \Sigma_{y,k} y^k\big)_i \in [(b_l)_i, (b_u)_i], \\
        \big(\Sigma_{y,k}\big)_{ii} \hspace{-0.2em} & \text{otherwise.}
    \end{cases}
\end{equation}
Therefore, the semismooth Newton direction $(\Delta x, \Delta \lambda)$ for
inner problem \eqref{eq:inner-qp} is the solution to
\begin{equation} \label{eq:semismooth-system}
    \begin{pmatrix}
        H_k(x) & \tp M \\ M & 0
    \end{pmatrix}
    \begin{pmatrix}
        \Delta x \\
        \Delta \lambda
    \end{pmatrix} = -\begin{pmatrix}
        \grad\phi_k(x) + \tp M \lambda \\
        Mx - b
    \end{pmatrix}.
\end{equation}

\subsection{Solving the semismooth Newton system}\label{subsec:semismooth-newton}

We use the approach described in \cite[Sec.~IV]{10550920} to solve
\eqref{eq:semismooth-system} in a way that preserves the block structure of the
original OCP matrices $Q, M$ and $G$:
\begin{equation} \label{eq:ocp-newton-system}
    \begin{aligned}
        H_k(x) v &= \grad\phi_k(x) + \tp M \lambda, \\
        -M \inv H_k(x) \tp M \Delta \lambda &= Mv - (Mx - b), \\
        -H_k(x) \Delta x &= \grad\phi_k(x) + \tp M (\lambda + \Delta \lambda).
    \end{aligned}
\end{equation}
Thanks to the structure of the matrices in \eqref{eq:standard-qp}, the matrix
$H_k(x)$ is block-diagonal, with blocks of the form
$H_{k,j}(x) = \begin{pmatrix}
    Q_j & \tp S_j \\ S_j & R_j
\end{pmatrix} + \begin{pmatrix}
    \tp C_j \\ \tp D_j
\end{pmatrix} \Sigma_{y,k,j}^{\setJ\!(x)} \begin{pmatrix}
    C_j & D_j
\end{pmatrix}$. This structure allows the evaluation and factorization of
$H_{k,j}(x)$ to be carried out in parallel across stages. The Cholesky
factorization $L_H \tp L_H$ of $H_{k,j}(x)$ can then be used to evaluate the
intermediate matrices $V_j$ and $W_j$ that produce the matrix
$\Psi_k(x) \defeq M \inv H_k(x) \tp M$. $\Psi_k(x)$ is block-tridiagonal,
with block-bidiagonal Cholesky factors consisting of $N+1$ diagonal blocks of
size $n_x\times n_x$, and $N$ coupling
blocks between stages along the subdiagonal. For the sake of readability, we
will from now on omit the dependence on $x$ and the iteration number $k$.
The following equation lists the
stage-wise parallelizable operations and the corresponding BLAS/LAPACK operations:
\begin{equation*}
    \begin{aligned}
        H_{j} &= \begin{pmatrix}
            Q_j & \tp S_j \\ S_j & R_j
        \end{pmatrix} + \begin{pmatrix}
            \tp C_j \\ \tp D_j
        \end{pmatrix} \Sigma_{y,j}^{\setJ} \begin{pmatrix}
            C_j & D_j
        \end{pmatrix} \hspace{-7em} & \texttt{(syrk)}\\
        L_{H,j} \tp L_{H,j} &= H_{j} & \texttt{(potrf)}\\
        V_j &= \begin{pmatrix} A_j & B_j \end{pmatrix} L_{H,j}^{-\top} & \texttt{(trsm)}\\
        W_j &= \begin{pmatrix} \I_{n_x} & 0_{n_x \times n_u} \end{pmatrix} L_{H,j}^{-\top} & \texttt{(trtri,gemm,trsm)}\\
        \Psi_{0,0} &= W_{0}\tp W_{0} & \texttt{(syrk)}\\
        \Psi_{j+1,j+1} &= V_j \tp V_j + W_{j+1}\tp W_{j+1} & \texttt{(syrk)}\\
        \Psi_{j+1,j} &= -V_j \tp W_{j} & \texttt{(gemm)}\\
    \end{aligned}
\end{equation*}
Next, the factorization of $\Psi$ is carried out recursively:
\begin{equation*}
    \begin{aligned}
        L_{\Psi,0,0}\tp L_{\Psi,0,0} &= \Psi_{0,0} & \texttt{(potrf)}\\
        L_{\Psi,j+1,j} &= \Psi_{j+1,j} L_{\Psi,j,j}^{-\top} & \texttt{(trsm)}\\
        \tilde \Psi_{j+1,j+1} &= \Psi_{j+1,j+1} - L_{\Psi,j+1,j}\tp L_{\Psi,j+1,j} \hspace{1.25em} & \texttt{(syrk)}\\
        L_{\Psi,j+1,j+1}&\tp L_{\Psi,j+1,j+1} = \tilde \Psi_{j+1,j+1} & \texttt{(potrf)}\\
    \end{aligned}
\end{equation*}
Finally, the vectors $v$, $\Delta \lambda$ and $\Delta x$ from \eqref{eq:ocp-newton-system}
can be computed using forward and back substitution of $L_{H}$ and $L_\Psi$.

All computations up to the factorization of $\Psi$ are fully independent and
will be parallelized in the following sections. The factorization of $\Psi$
can still make use of the parallel nature of linear algebra operations inside
of each $n_x \times n_x$ block.

\section{Vectorization} \label{sec:vec}

Many of the operations in QPALM-OCP can be expressed stage-wise. This includes
part of the solution of \eqref{eq:ocp-newton-system} as discussed in the
previous section, but also operations such as the matrix-vector product $Mx$,
where one needs to evaluate the stage-wise products $A_j x^j$. These smaller
products are fully independent and
can be carried out in parallel. Since the same operation is applied to different
matrices, we can use single instruction, multiple data (SIMD) programming
techniques across the multiple stages.

\subsection{Compact storage format}
\begin{figure}
    \centering
    \includegraphics[width=\linewidth]{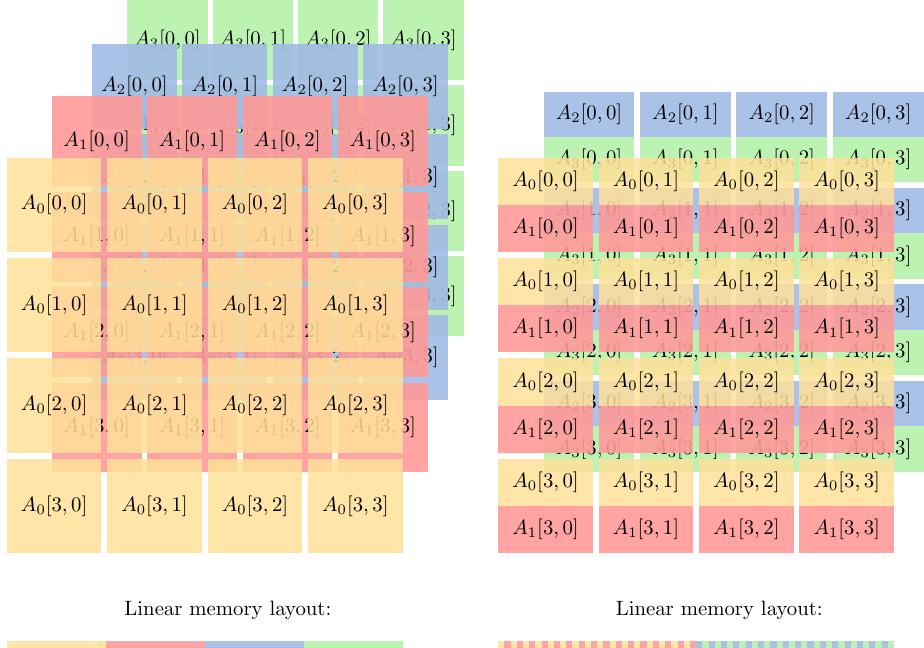}
    \caption{Comparison between ``naive'' (left) and ``compact'' storage (right) of matrices $A_j$ for a problem of size $N = 4 = n_x$ and a vector length of $d=2$.
    In the naive format, each matrix $A_j$ is stored contiguously in memory (e.g. column major order), and matrix $A_{j+1}$ of the next stage is stored after $A_j$.
    In the compact format, matrices $A_{2k}$ and $A_{2k+1}$ are interleaved in such a way that their elements $(A_{2k})_{\imath\jmath}$ and $(A_{2k+1})_{\imath\jmath}$ are stored next to each other.
    Alternatively, the naive format can be seen as a $n_x \times n_x \times N$ tensor where the stage number $j \lt N$ is the index with the largest stride. In contrast, the
    compact format would be a $n_x \times n_x \times \lceil N/d \rceil$ tensor where each element is a tuple of $d$ elements.} \label{fig:simd-format}\vspace{-1em}
\end{figure}

To enable this data parallelism across stages, the full OCP horizon $N$ is
subdivided into blocks with the same size as the hardware's vector length
(2, 4 or 8 for modern CPUs), and all stages within a block are processed simultaneously.
For example, using a vector length of 2, the products $A_0 x^0$ and $A_1 x^1$
are carried out at the same time, and so are the products $A_2 x^2$ and $A_3 x^3$.
In order to perform such SIMD operations efficiently, the corresponding elements
$(A_0)_{\imath\jmath}$ and $(A_1)_{\imath\jmath}$ need to be located at adjacent memory addresses.
To this end, we store the matrices in a ``compact'' format, where matrices of
different stages are interleaved. \Cref{fig:simd-format} visualizes this format.

In our implementation of QPALM-OCP, all OCP matrices ($A_j$, $B_j$, $C_j$,
$Q_j$, $R_j$, $S_j$) are stored in a compact format, and so are the
corresponding vectors ($x$, $\Sigma_y$, $y$, etc.), in order to minimize the
number of storage conversions required over the course of the algorithm.

When the horizon length $N$ is not a multiple of the vector length, zero or
identity matrices are added as appropriate, and the final stage is padded with
$S^N = 0, R^N = \I, D^N=0$.

\subsection{Vectorized linear algebra routines for matrices in compact storage}
Even though most linear algebra libraries will readily use SIMD instructions for
ordinary matrices, using SIMD \textit{across} multiple matrices rather than \textit{within} a single matrix
can result in higher performance, especially for smaller matrices \cite{intel_compact_perf}.
For this reason, a selection of linear algebra routines that operate on batches of
matrices in compact storage format are available in the Intel Math Kernel Library (MKL) \cite{mkl_compact}.
In our benchmarks, these MKL routines resulted in some avoidable overhead,
so we implemented some of these performance-critical routines
ourselves. The architecture is roughly based on the BLIS framework \cite{blis}
\cite[Fig.~1]{cgemm}, with
multiple layers of loops around optimized GEMM and TRSM-GEMM micro-kernels.
These micro-kernels keep a small $3\times 3\times d$ or $5\times 5\times d$ block of the
result matrices in the CPU's vector registers as an accumulator, and perform some
highly optimized multiply-accumulate operations. The loops
around the micro-kernel perform blocking (and optionally packing) to optimize
cache usage. 
The implementation is
built on top of the experimental
\texttt{std:\!:simd} library \cite{stdsimd} and makes use of \Cpp{} function templates for handling different block sizes and data types.
The practical effectiveness of the storage format and the optimized linear
algebra routines will be discussed in \Cref{sec:results}.

\section{Parallelization} \label{sec:par}
It is unlikely that the full horizon length $N$ fits inside of a single vector
register. Therefore we can further parallelize the $\lceil N/d\rceil$ blocks of
stages introduced in the previous section by distributing them across the
different CPUs in modern multi-core hardware. Frameworks like OpenMP \cite{openmp}
provide an easy way to express such parallelism.

\section{Numerical examples}\label{sec:results}
This section reports the performance of a \Cpp{} implementation of QPALM-OCP
using the optimizations described in the previous sections. \citet{10550920}
used the classical spring-mass benchmark from \cite{springmass} to compare the
performance of their prototype implementation of QPALM-OCP to QPALM, OSQP
and HPIPM. We keep QPALM as a baseline, and also include the
recent PIQP solver \cite{piqp}.

Furthermore, we also investigate the effectiveness of our parallelization efforts, and conclude
with four test problems from the \texttt{qpsolvers/mpc\_qpbenchmark} repository \cite{qpbenchmark2024}.

\subsection{Setup}
All benchmarks were carried out using an octa-core Intel Core i7-11700 @
2.50\,GHz with dynamic CPU frequency scaling disabled for deterministic results. Native code was compiled 
using GCC 14.1 with AVX-512 support enabled at optimization level \texttt{-O3}.
The absolute tolerances for all solvers are set to $10^{-8}$, and the relative
tolerances are set to zero.

\subsection{Spring-mass benchmark}
We use the same setup as described in \cite[Sec.~V]{10550920}: The number of
masses $M$ is varied from 10 to 70, and the horizon is fixed at $N=15$. For
each configuration, we solve 100 OCPs with randomly generated initial conditions
as in \cite{10550920}, and the geometric means of the solver run times are reported.

Since this specific problem has diagonal cost matrices $Q_j$ and $R_j$, as well as
diagonal constraint matrices $\begin{pmatrix} C_j \; D_j \end{pmatrix}$ (box constraints),
the resulting matrices $H_{k,j}(x)$ are diagonal as well, significantly simplifying
the solution of \eqref{eq:ocp-newton-system}. We therefore compare two different
implementations of the method described in \Cref{subsec:semismooth-newton}:
One general implementation where all problem data is represented by full, dense matrices,
and one implementation where we exploit the diagonal structure of $H_{k,j}$ as much as possible.

For a fair comparison, we also include two versions of the sparse QPALM and PIQP
solvers: one where we include all OCP matrices as dense blocks in the sparse
matrix (even though they contain numerical zeros), and one where we eliminate
all zeros from the sparse matrices. In the case of PIQP, we also provide the
constraints as bounds on the variables rather than using general linear
constraints for the second scenario.

\begin{figure}[t]
    \centering
    \includegraphics[width=\linewidth]{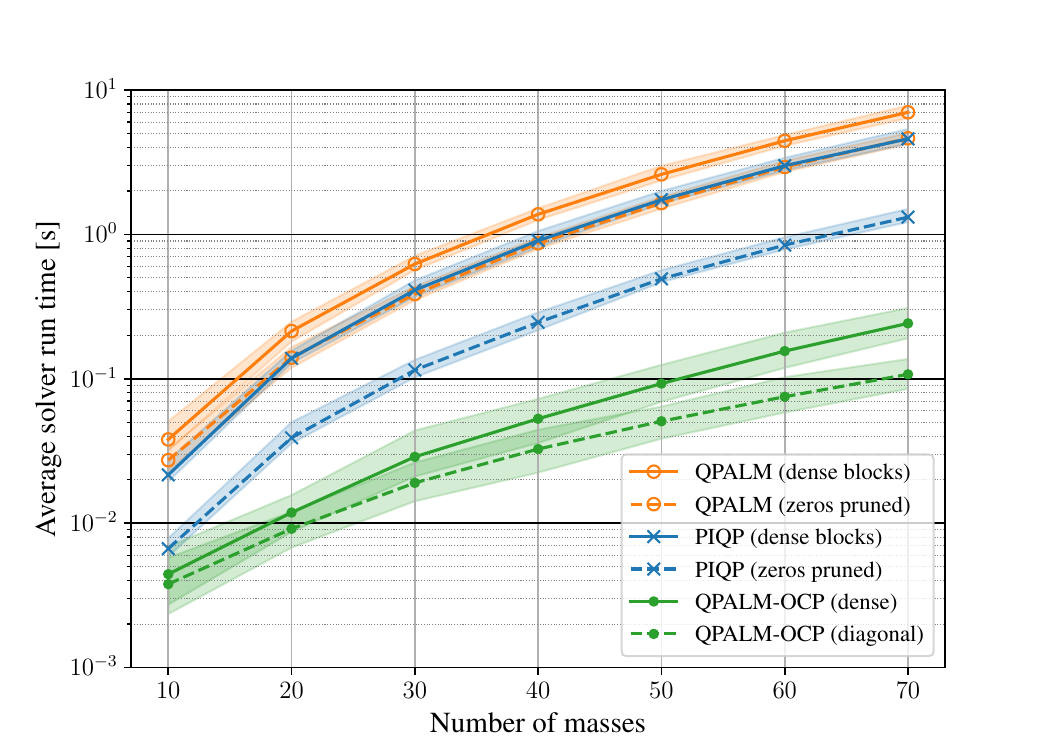}
    \vspace{-1em}
    \caption{Average solver run times for different numbers of masses, horizon $N=15$. Shaded areas indicate best and worst run times.} \label{fig:spring-masses}
\end{figure}

The results of this experiment are shown in \Cref{fig:spring-masses}.
Both optimized versions of QPALM-OCP outperform the other solvers on this benchmark. For the largest problem (3275 primal variables), the dense version of QPALM-OCP is around $29$ times faster than QPALM with dense blocks, and over $19$ times faster than QPALM with all numerical zeros pruned.
The diagonal version of QPALM-OCP is around $65$ times faster than QPALM with dense blocks, and over $43$ times faster than QPALM with all numerical zeros pruned.

\subsection{Effectiveness of parallelization}
\begin{figure}[t]
    \centering
    \includegraphics[width=\linewidth]{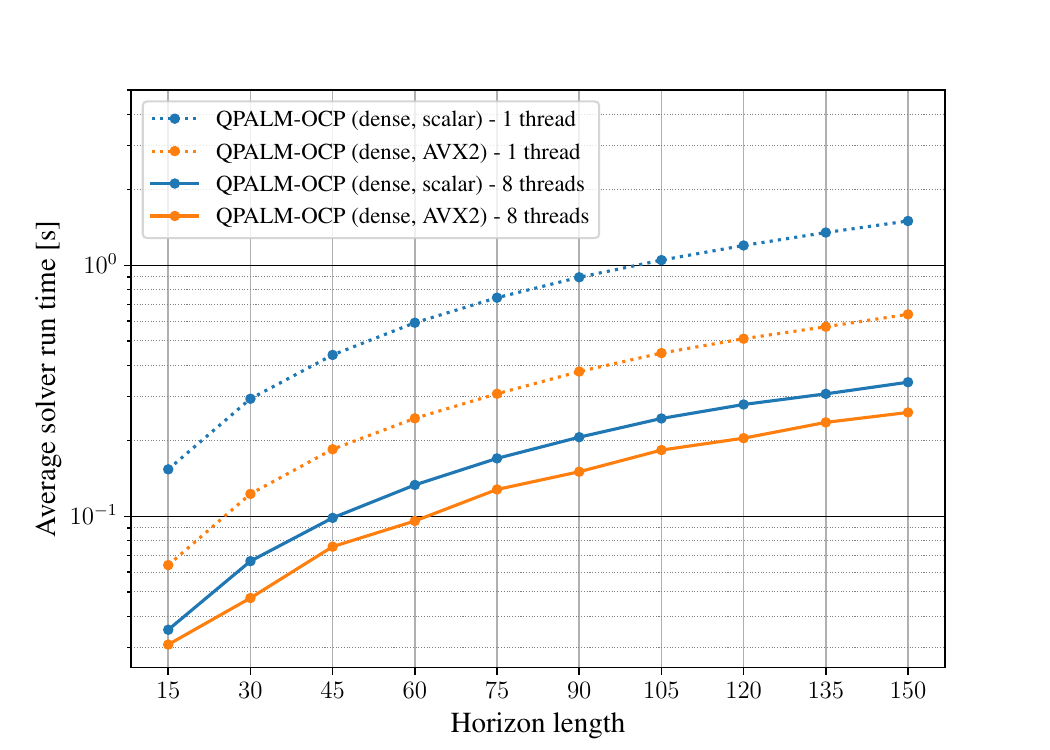}
    \vspace{-1em}
    \caption{Effect of parallelization and vectorization for different horizon length, masses $M=30$.} \label{fig:spring-masses-simd}
\end{figure}

In this experiment, we fix the number of masses to $M=30$ and run four variants of the solver for different horizon lengths. The dense QPALM-OCP solver is used with both AVX2 vectorization (vector length 4), and without vectorization (no compact storage). These two configurations are repeated for a single OpenMP thread, and 8 OpenMP threads. The geometric means of the run times for 10 random instances of the spring-mass benchmark are reported.
\Cref{fig:spring-masses-simd} summarizes the results:
In the single-threaded case, vectorization results in an overall speedup by a factor of around 2.3. By using multiple threads, the performance increases further, but it is limited by cache and bandwidth effects as well as the sequential parts such as the factorization of $\Psi$.

\subsection{MPC qpbenchmark}
Finally, we briefly report the performance on the four \texttt{QUADCMPC*} quadruped locomotion benchmarks from the \texttt{qpsolvers/mpc\_qpbenchmark} repository \cite{qpbenchmark2024} in \Cref{tab:performance}. Even though the problem is very sparse (yet not diagonal), the dense QPALM-OCP solver significantly outperforms the sparse QPALM solver.
On the \texttt{LIPMWALK*} benchmarks, the performance difference is less pronounced, since these problems are very small ($n_x=3, n_u=1, N=16$). With OpenMP disabled (it adds overhead for these small problems), QPALM-OCP achieves a shifted geometric mean run time of 0.43\,ms, versus 0.46\,ms for standard QPALM.

\begin{table}[h]
    \centering
    \begin{tabular}{lrr}
    \toprule
    Problem & QPALM (pruned) & QPALM-OCP (dense) \\
    \midrule
    \texttt{QUADCMPC1} & 21.2\,ms & 5.1\,ms \\
    \texttt{QUADCMPC2} & 9.9\,ms & 2.7\,ms \\
    \texttt{QUADCMPC3} & 3.0\,ms & 2.0\,ms \\
    \texttt{QUADCMPC4} & 3.5\,ms & 0.7\,ms \\
    \bottomrule
    \end{tabular}
    \caption{\normalfont Performance comparison between QPALM and QPALM-OCP on the \texttt{qpsolvers/mpc\_qpbenchmark} \texttt{QUADCMPC*} problems.}
    \label{tab:performance}
    \vspace{-2em}
\end{table}

\section{Conclusion} 
\label{sec:conclusion}

In this paper, we have shown that the recent QPALM-OCP method can effectively be parallelized, and that it outperforms other state-of-the-art solvers on a set of benchmark problems.
Future work includes efficient offline packing of the matrix storage \cite{blasfeo,blis}, and the implementation of factorization update routines to avoid complete refactorizations when a small number of constraints or penalty factors change.

\clearpage
\bibliographystyle{plainnat}
\bibliography{references}

\section*{Acknowledgments}

This work was supported by Fonds Wetenschappelijk Onderzoek (FWO)
PhD grant 11M9523N; FWO research projects G081222N, G033822N and
G0A0920N; the European Union's Horizon 2020 research and
innovation programme under the Marie Sk\l{}odowska-Curie
grant agreement No. 953348; and EuroHPC Project 101118139 Inno4Scale.

\end{document}

%% file: mynotation.tex
\newcommand\defeq{\triangleq}
\newcommand\iddots{\mathinner{
    \kern1mu\raise1pt{.}
    \kern2mu\raise4pt{.}
    \kern2mu\raise7pt{\Rule{0pt}{7pt}{0pt}.}
    \kern1mu
}}

\newcommand\lt{<}

\newcommand\grad{\nabla}

\newcommand\Reals{\mathrm{I\!R}}

\newcommand\R{\Reals}

\newcommand\possdefset[1]{\operatorname{Sym}_+(#1)}
\newcommand\posdefset[1]{\operatorname{Sym}_{++}(#1)}

\newcommand\minimize{\operatorname*{\mathbf{minimize}}}

\newcommand\subjto{\operatorname*{\mathbf{subject\;to}}}

\newcommand\dist{\operatorname{\mathbf{dist}}}

\newcommand\I{\mathrm{I}}

\newcommand\norm[1]{\left\| {#1} \vphantom{X} \right\|}

\newcommand\normsq[1]{\norm{#1}^2}

\newcommand\defset[2]{\left\{ {#1} \;\middle|\; {#2} \right\}}
\newcommand\inv[1]{#1^{-1}}

\newcommand\tp[1]{#1^\top}

\newcommand\Cpp{C\texttt{++}}